\documentclass[11pt]{amsart}
\usepackage{epsfig}

\newcommand{\lcr}{\raisebox{-5pt}{\mbox{}\hspace{1pt}
                  \epsfig{file=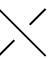}\hspace{1pt}\mbox{}}}
\newcommand{\ift}{\raisebox{-5pt}{\mbox{}\hspace{1pt}
                  \epsfig{file=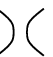}\hspace{1pt}\mbox{}}}
\newcommand{\zer}{\raisebox{-5pt}{\mbox{}\hspace{1pt}
                  \epsfig{file=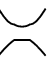}\hspace{1pt}\mbox{}}}

\title[The A-polynomial and the Jones polynomial]{On the relation between 
the A-polynomial and the Jones polynomial}

\author{R{\u{a}}zvan Gelca}
\address{Department of Mathematics,University of Michigan, Ann Arbor, MI 48109, USA and Institute of Mathematics of the Romanian Academy, 
Bucharest, Romania}
\email{rgelca@math.lsa.umich.edu}

\newtheorem{theorem}{Theorem}

\newtheorem{lemma}{Lemma}

\begin{document}
\maketitle
\section{Introduction}

In 1984, V. Jones introduced a polynomial invariant of knots
\cite{J} through skein relations. Another version
of this invariant, the Kauffman bracket \cite{K}, was introduced shortly
 after. Colored versions of these invariants were defined,
 via quantum groups \cite{RT}, and via Jones-Wenzl idempotents \cite{L},
\cite{BHMV}.

In 1993, Cooper, Culler, Gillet, Long, and Shalen defined a two variable 
polynomial invariant of knots, the A-polynomial,
 using the character variety of
$SL(2,{\mathbb C})$-representations of the fundamental group of the knot
complement. This invariant was generalized in \cite{FGL} to 
a finitely generated ideal of polynomials in the quantum plane.
The construction is done in the context of skein modules,
and is based on the fact that the  Kauffman bracket 
skein modules represent deformations of function rings on 
character varieties \cite{B}, \cite{SP}, and on the relationship
between the skein algebra of the cylinder over a torus and
the noncommutative torus \cite{FG}. 

As shown in \cite{FGL}, each element in the noncommutative A-ideal defines
a matrix that annihilates the vector whose entries are the colored Jones
polynomials of the knot (or, more precisely, the colored Kauffman brackets
of the knot; they differ from the colored Jones polynomials by the
change of variable $t\mapsto it$).
 The orthogonality between the rows of the matrix
and the vector whose entries are the colored Jones polynomials of the knot
has been called the ``orthogonality relation (between the Jones polynomial
and the A-polynomial)''. 

In the present paper it is shown that the noncommutative A-ideal together
with a finite number (depending on the A-ideal of the  knot)
 of colored Kauffman brackets  of the knot determine all   
other colored Kauffman brackets of the knot. 
Also, it is shown that, under certain technical conditions on the
A-ideal, the noncommutative A-ideal determines all colored Kauffman
brackets of the knot. As an example, any knot having the same
A-ideal as the unknot, respectively trefoil knot, has the same
colored Kauffman brackets as the unknot, respectively trefoil knot.

\section{The action of $K_t({\mathbb T}^2\times I)$ on
$K_t({\mathbb D}^2\times I)$}

The Kauffman bracket skein module of the 
three manifold $M$ is defined in the following way.
 Let ${\mathbb C}[t,t^{-1}]{\mathcal L}$ be the  
${\mathbb C}[t,t^{-1}]$-module  freely spanned by the 
isotopy classes of framed links in $M$ including the 
empty link, and let ${\mathcal S}$ be the submodule 
spanned by the  relations $\displaystyle{\lcr-t\zer-t^{-1}\ift}$
and 
$\bigcirc+t^2+t^{-2}$. The Kauffman
bracket skein module  of $M$ is $K_t(M)={\mathbb C}[t]/
{\mathcal S}$.

 In the case where $M$ is the cylinder over 
a surface, $K_t(M)$ has a natural algebra structure, with
product defined by placing one link on top of another.
If $M$ is a manifold with boundary, the 
operation of gluing a cylinder to the boundary 
 induces a $K_t(\partial M\times I)$-module 
structure on $K_t(M)$.
As an  example it is known that the Kauffman bracket
skein  algebra of the cylinder
over an annulus (i.e., that of the solid torus),
 is ${\mathbb C}[t,t^{-1},\alpha]$,
where $\alpha $ is the curve that runs once around the annulus
and has framing parallel to the annulus.

Another, more complicated example is that of the 
Kauffman bracket skein algebra of $K_t({\mathbb T^2}\times I)$.
Its multiplication rule
 and action on the skein module of the 
solid torus are described by means of two families of
Chebyshev polynomials, $\{T_n\}_{n\in {\mathbb Z}}$ defined by $
T_0=2, T_1=x,$ $T_{n+1}=xT_n-T_{n-1}$ for $n\in {\mathbb Z }$
and $\{S_n\}_{n\in {\mathbb Z}}$ defined by $
S_0=2, S_1=x,$ $S_{n+1}=xS_n-S_{n-1}$ for $n\in {\mathbb Z }$.
Let $p$ and $q$ be two integers with $p=np'$, $q=nq'$,
$p'$, $q'$ coprime. We define $(p,q)_T=T_n((p',q'))$, where
$(p',q')$ is the corresponding curve on the torus, with framing
parallel to the torus, and its powers are defined by parallel copies.
The elements $(p,q)_T$, $p\geq 0, q\in {\mathbb Z}$ span
$K_t(M)$ as a ${\mathbb C}[t,t^{-1}]$-module. 
In \cite{FG} we proved  the following product-to-sum formula
\[
(p,q)_T*(r,s)_T=t^{|^{pq}_{rs}|}
(p+r, q+s)_T+
t^{-|^{pq}_{rs}|}
(p-r,q-s)_T.
\]

As a consequence of this formula, the Kauffman bracket 
skein algebra of the cylinder over a torus is isomorphic 
to the subalgebra of the noncommutative torus generated
by noncommutative cosines. Let us recall that the algebra of 
trigonometric polynomials in the noncommutative torus is
$ {\mathbb C}_t[l,l^{-1}, m,m^{-1}]$, with multiplication 
$*$, satisfying  $l*m=t^2m*l$. 
The elements $e_{p,q}=t^{-pq}l^pm^q$ are the noncommutative 
exponentials;
they satisfy 
\begin{eqnarray*}
e_{p,q}*e_{r,s}=t^{|^{pq}_{rs}|}
e_{p+r,q+s}.
\end{eqnarray*}
The noncommutative  cosines are $\frac{1}{2}(e_{p,q}+e_{-p,-q})$.
The map $(p,q)_T\rightarrow e_{p,q}+e_{-p,-q}$ gives the 
isomorphism between $K_t({\mathbb T}^2\times I)$ and
the algebra of noncommutative cosines.        

Let $K$ be a knot in $S^3$, and $M$ the complement of a regular
neighborhood of $K$. Recall the left action of 
$K_t({\mathbb T}^2\times I)$ on $K_t(M)$. 
The peripheral ideal of $K$ is the 
left ideal of $K_t({\mathbb T}^2\times I)$ which 
annihilates the empty link. The noncommutative A-ideal of 
$K$, denoted by ${\mathcal A}_t(K)$
 is the left ideal obtained by extending $I_t(K)$ to
${\mathbb C}_t[l,l^{-1},m,m^{-1}]$ then contracting
it to  ${\mathbb C}_t[l, m]$. As explained
in \cite{FGL}, this is a noncommutative generalization of
the A-polynomial. The A-polynomial is obtained by setting 
$t=-1$,  replacing $l$ and $m$ by $-l$ and $-m$ and taking
the generator of the radical of the one-dimensional part of the
A-ideal (divided by $(l-1)$).

There is a left  and a right action of $K_t({\mathbb T}^2\times
I)$ on $K_t({\mathbb D}^2\times S^1)$, one for the 
positive, the other one for the negative orientation 
of the boundary torus. To understand them, let us denote
by $x_{p,q}$ the image  in $K_t({\mathbb D}^2\times S^1)$
of $(p,q)_T$ on the boundary torus (with the positive orientation).
It is not hard to see that $x_{0,q}= (-t^2)^q+(-t^{-2})^q$ and
the product-to-sum formula yields
\begin{eqnarray*}
x_{p+1,q}=t^{-q}(1,0)\cdot x_{p,q}-t^{-2q}x_{p-1,q}.
\end{eqnarray*}
The second order recurrence relation for $t^{pq}x_{p,q}$ has
fixed coefficients, and hence a formula for  the general term
can be found. It is 
\begin{eqnarray*}
x_{p,q}=t^{-pq}((-t^{-2})^qS_p(\alpha)-
(-t^2)^qS_{p-2}(\alpha)).
\end{eqnarray*}

Lifting the skeins $T_n(\alpha )$ to the boundary torus
and using the product-to-sum formula we get the following 
\begin{lemma}
The left action is described by
\begin{eqnarray*}
(p,q)_T\cdot T_n(\alpha)& = &  t^{-(2n+p)q}[(-t^{-2})^q
S_{n+p}(\alpha )-(-t^{-2})^qS_{n+p-2}(\alpha )]\\ & & +t^{(2n-p)q}[(-t^{-2})^q
S_{p-n}(\alpha )-(-t^{2})^qS_{p-n-2}(\alpha)]
\end{eqnarray*}
while the right action is given by
\begin{eqnarray*}
 T_n(\alpha )\cdot (p,q)_T & = & (p,-q)_T\cdot T_n(\alpha )\\ & = & 
t^{(2n+p)q}[(-t^2)^qS_{p+n}(\alpha)
 -(-t^{-2})^qS_{p+n-2}(\alpha)]\\ & & 
+t^{-(2n-p)q}[(-t^2)^qS_{p-n}(\alpha)-
(-t^{-2})^qS_{p-n-2}(\alpha )].
\end{eqnarray*}
\end{lemma}

\section{The results}

Gluing a solid torus to the  complement $M$ of a regular neighborhood
of a knot $K$, in such
a way that the longitude is glued to the longitude and
the meridian to the meridian, induces a pairing
\begin{eqnarray*}
K_t({\mathbb D}\times S^1)\times K_t(M)\rightarrow {\mathbb C}
[t,t^{-1}].
\end{eqnarray*}
 The basis $\{S_n(\alpha )\}_n$ induces a family of functionals
$< S_n(\alpha ), \cdot >$, $n=0,1,2,\ldots $. If we denote by 
$\emptyset$ the empty link, then
\begin{eqnarray*}
<S_n(\alpha ), \emptyset >=\kappa _n(K),
\end{eqnarray*}
where $\kappa_n(K)$ is the $n$th colored Kauffman bracket 
of $K$ with zero framing \cite{L}, \cite{T} (the $n$th colored Kauffman 
bracket is a
``twisted'' version of the $n$th  colored Jones
 polynomial as defined in \cite{RT}). 
Indeed, the recurrence relation
for $S_n$ shows that the link in $S^3$  obtained from the
pairing  is $K$ colored by the Jones-Wenzl idempotent.

The pairing is compatible with the actions of 
$K_t({\mathbb T}^2\times I)$ on both modules, i.e.
$<u\cdot (p,q)_T,v>=<u, (p,q)_T\cdot v>$ for any
skeins $u$ and $v$. 
In particular, if $a$ is in the peripheral ideal $I_t(K)$
of $K$, then $<u\cdot a, \emptyset >=0$. So,
if $u=T_n(\alpha)$, and $a=\sum_i c_i(p_i,q_i)_T$, then
by Lemma 1,
\begin{eqnarray*}
& & <a\cdot T_n(\alpha ),\emptyset>  =   \sum_ic_i(t^{p_iq_i}
(t^{(2n+p_i)q_i}[(-t^2)^{q_i}<S_{p_i+n}(\alpha),\emptyset>\\ & &  -
(-t^{-2})^{q_i}<S_{p_i+n-2}(\alpha),\emptyset>]\\ & & 
+t^{-(2n-p_i)q_i}[(-t^2)^{q_i}<S_{p_i-n}(\alpha), \emptyset> -
(-t^{-2})^{q_i}<S_{p_i-n-2}(\alpha ),\emptyset>])\\
 & = & \sum_ic_i(t^{(2n+p_i)q_i}[(-t^2)^{q_i}\kappa_{p_i+n}(K)
 -(-t^{-2})^{q_i}\kappa_{p_i+n-2}(K)]\\ & & +t^{-(2n-p_i)q_i
}[(-t^2)^{q_i}\kappa_{p_i-n}(K)-
(-t^{-2})^{q_i}\kappa_{p_i-n-2}(K)]).
\end{eqnarray*}
This relation has been called the orthogonality relation in
\cite{FGL} since
it expresses the orthogonality between the vector with
entries equal to the colored Kauffman brackets of the knot and
the rows of the matrix of the linear transformation 
induced by $a$ between the module $K_t({\mathbb D}\times I)$ with
basis $\{T_n(\alpha)\}_n$ and the same module with basis
$\{S_n(\alpha)\}_n$. Since $a$ arises from an element in
the noncommutative A-ideal (through an extension and a contraction),
 orthogonality expresses a relationship between the 
the elements of the A-ideal and the vector whose entries are
the colored Kauffman brackets.  
\begin{theorem}
For every knot $K$ there is a number $\nu (K)$ such that if 
$K'$ is a knot with ${\mathcal A}_t (K)={\mathcal A}_t(K')$ and
$\kappa _j(K)=\kappa _j(K')$ for $j=1,2,\cdots , \nu(K)$, then
$\kappa _j(K)=\kappa _j(K')$ for all $j$. Moreover, $\nu(K)$ 
depends only on the A-ideal of $K$.
\end{theorem}

\proof Choose $a=\sum_j c_{j}(p_j,q_j)_T$ some element in 
 $I_t (K)$, let $p$ be the maximum of $p_j$ and assume
$p_j=p $ if $j=1,2,\cdots , m$, $p_j\neq p$ if $j>m$. Then,
 coefficient of $\kappa _{n+p}$ in the orthogonality relation
written for $a$ is 
\begin{eqnarray*}
\sum_{j=1}^mc_j(-1)^{q_j}t^{(2n+2+p)q_j}
\end{eqnarray*}
Since the $q_j$ appearing in this expression are distinct (the
$p_i$'s being the same), this expression is identically equal to
zero only for finitely many $n$. Hence the orthogonality relation
provides a recurrence relation that determines uniquely 
$\kappa _n$ for large $n$. \qed

As the result below shows, in certain situations the A-ideal 
determines the colored Kauffman brackets of the knot.

\begin{theorem}
Assume that $K$ is a knot with the property that ${\mathcal A}_t(K)$
contains a polynomial $\sum _{p,q}\gamma _{p,q}l^pm^q$
 of degree $2$ in $l$ such that there exists no $n\geq 0$ for which the 
expression 
$\sum_{q}\gamma _{2,q}(-1)^qt^{(2n+2)q}$ is identically equal to zero.
 Then for any knot 
$K'$ with the property that ${\mathcal A}_t(K)={\mathcal A}_t(K')$,
it follows that $\kappa_n(K)=\kappa_n(K')$ for all $n=1,2,3,\ldots$.
\end{theorem}

\proof The  polynomial gives rise to an element $a=\sum_ic_i(1,q_i)_T+u$
in $I_t(K)$, with $c_i=t^{q_i}\gamma _{2,q_i}$ and 
$u$ a polynomial in $(0,1)$.
By Lemma 1, $T_n(\alpha )\cdot u$ is of the form
$\lambda S_n(\alpha )+\mu S_{n-2}(\alpha)$, $\lambda ,\mu \in
{\mathbb C}[t,t^{-1}]$. On the other hand, the same lemma shows that
\begin{eqnarray*}
& & T_n(\alpha)\cdot \sum_ic_i(1,q_i)_T=\sum_i c_i[(-t)^{(2n+3)q_i}
S_{n+1}(\alpha )\\ & & -(-t)^{(2n-1)q_i}S_{n-1}(\alpha)+(-t)^{(-2n+3)}
S_{1-n}(\alpha)-(-t)^{(-2n-1)q_i}S_{-n-1}(\alpha)].
\end{eqnarray*}
Since $S_{-k}=-S_{k-2}$ for all $k$, 
this is further equal to
\begin{eqnarray*}
 \sum _ic_i[(-t)^{(2n+3)q_i}S_{n+1}(\alpha )-
[(-t)^{(2n-1)q_i}+(-t)^{(-2n-1)q_i }]S_{n-1}(\alpha )\\
+(-t)^{(-2n+3)q_i}S_{n-3}
(\alpha)].
\end{eqnarray*}
Hence the orthogonality relation applied to $a$ yields a
${\mathbb C}[t,t^{-1}]$-linear equation in
$\kappa _{n+1}(K), \kappa _n(K), \kappa_{n-1}(K), \kappa _{n-2}(K)$,
 and $\kappa _{n-3}(K)$.
 For $n\geq 1$, the coefficient of $\kappa _{n+1}(K)$ is
$\sum_ic_i(-t)^{(2n+3)q_i}$, and the condition from the statement translates
to the fact that for no $n$ this is identically equal to zero.
 Therefore,  the orthogonality 
relation provides a recursion that uniquely determines $\kappa _n(K)$ from
$\kappa _1(K)$. 
 
For $n=0$, since $\kappa _0(K)=1$, 
$\kappa _{-1}(K)=0$, $\kappa _{-2}(K)=-1$, $\kappa _{-3}(K)=
-\kappa _1(K)$, 
the orthogonality relation gives a linear equation in 
$\kappa _1(K)$. The coefficient of $\kappa _1(K)$ is 
$2\sum_ic_i(-t)^{3q_i}$. Again this is not equal to zero.
So the equation  
can be  solved uniquely for $\kappa _1(K)$.
It follows that the A-ideal determines the colored Kauffman brackets
of the knot, and we are done.\qed

Observe that the degree in $l$ 
of any polynomial in ${\mathcal A}_t(K)$ is at least
$2$.

\section{Examples}

\subsection{The unknot}

The A-ideal of the unknot is generated by 
$(l+t^2)(l+t^{-2})$ and $lm^2(l+t^2)+t^2(l+t^{-2})$  \cite{FGL};
hence it satisfies the conditions in Theorem 2. The orthogonality 
relation for $(l+t^2)(l+t^{-2})$, that is, for $(1,0)_T+t^2+t^{-2}\in I_t(K)$,
  gives  
\begin{eqnarray*}
& \kappa _{0}(K)=1, \quad \kappa _1(K)=-t^2-t^{-2}, & \\
&  \kappa_{n+1}
(K)=(-t^2-t^{-2})\kappa _n(K)-\kappa _{n-1}(K), \quad n\geq 1.&
\end{eqnarray*}
From this  we obtain the well known formula 
\begin{eqnarray*}
\kappa _n(K)=(-1)^n(t^{2n+2}-t^{-2n-2})/(t^2-t^{-2}).
\end{eqnarray*}
The orthogonality relation for the other element leads to a different
recurrence relation with the same solution.

\subsection{The trefoil}

The A-ideal of the left-handed trefoil is generated by 
$[m^4(l+t^{10})-t^{-4}(l+t^{2})](l-t^6m^6)$, 
$(l+t^{24})(l+t^{10})(l+t^2)(l-t^6m^6)$ and 
$(m^2-t^{-22})(l+t^{10})(l+t^2)(l-t^6m^6)$  \cite{G}.
A quick look at the element $[m^4(l+t^{10})-t^{-4}(l+t^{2})](lm^6-t^6)$
 shows that the conditions in the statement of Theorem 2 are fulfilled.
This element  corresponds to 
\begin{eqnarray*}
(1,-5)_T-t^{-8}(1,-1)_T+t^3(0,5)_T-t(0,1)_T
\end{eqnarray*}
in the peripheral ideal.
The orthogonality relation produces the following 
recursion 
\begin{eqnarray*}
& (-t^{-10n-15}+t^{-2n-11})\kappa _{n+1}(K)+
(-t^{10n+7}-t^{-10n-13}+t^{2n+3} & \\ & 
+t^{-2n-1})\kappa _n(K) +(t^{-10n+5}-t^{10n+5}-t^{-2n-7}+t^{2n-7})\kappa _{n-1}(K) +(t^{10n-13}& \\ &
+t^{-10n+7}-t^{2n-1}-t^{-2n+3})\kappa _{n-2}(K)+
(t^{10n-15}-t^{2n-11})\kappa _{n-3}(K)=0.&
\end{eqnarray*} 
In particular, for $n=0$,  
\begin{eqnarray*}
(t^{-11}-t^{-15})\kappa _1(K)-t^7-t^{-13}+t^3+t^{-1}=0,
\end{eqnarray*}
and hence $\kappa _1(K)=t^{18}-t^{10}-t^6-t^2$, the well known
formula for the Kauffman bracket of the trefoil knot with 
framing zero.

\end{document}